\documentclass[11pt]{amsart}

\usepackage{verbatim} 

\usepackage{fancyhdr}

\usepackage{amsmath, amsthm, amssymb}

\usepackage{graphicx}

\usepackage{amstext}

\def\Z{{\mathbb Z}}

\def\INT{{\mathrm{int}}}

\def\TSG{{\mathrm{TSG_+}}}
\def\Aut{{\mathrm{Aut}}}
\def\Diff{{\mathrm{Diff_+}}}
\def\fix{{\mathrm{fix}}}

\def\TSG{{\mathrm{TSG_+}}}
\def\Aut{{\mathrm{Aut}}}
\def\Diff{{\mathrm{Diff_+}}}
\def\fix{{\mathrm{fix}}}

\newtheorem*{complete}{Complete Graph Theorem}
\newtheorem*{finiteness}{Finiteness Lemma}
\newtheorem*{auto}{Automorphism Theorem}
\newtheorem*{orbits}{Orbits Lemma}
\newtheorem*{disjoint}{Disjoint Fixed Points Lemma}
\newtheorem*{edge}{Edge Embedding Lemma}

\newtheorem*{fixed}{Fixed Vertex Lemma}
\newtheorem*{3cycle}{3-Cycle Lemma}
\newtheorem*{pq}{pq Lemma}
\newtheorem*{4r+3}{ Theorem}
\newtheorem*{D2}{$D_2$ Lemma}

\def\fix{{\mathrm{fix}}}

\newtheorem{definition}{Definition}

\newtheorem*{subgroup}{Subgroup Theorem}

\title[Topological symmetry groups of $K_{4r+3}$]{Topological symmetry groups of $K_{4r+3}$}

\author{Dwayne Chambers}
\address{Department of Mathematics,
  Claremont Graduate University, Claremont CA 91711, USA.}
\email{dwayne.chambers@cgu.edu}
\author{Erica Flapan}
\address{Department of Mathematics, Pomona College, Claremont CA 91711, USA.}
\email{eflapan@pomona.edu}
\author{John D. O'Brien}
\address{Centre for Genomics and Global Health, Oxford University, Oxford OX3 7BN, UK.}
\email{jobrien@well.ox.ac.uk}



\thanks{The first two authors were partially supported by NSF grant DMS-0905087.}

\begin{document}

\maketitle

\begin{abstract}  We present the concept of the topological symmetry group as a way to analyze the symmetries of non-rigid molecules.  Then we characterize all of the groups which can occur as the topological symmetry group of an embedding of complete graphs of the form $K_{4r+3}$ in $S^3$.
\end{abstract}

\section{Topological Symmetry Groups}

\begin{figure} [htp]
\includegraphics{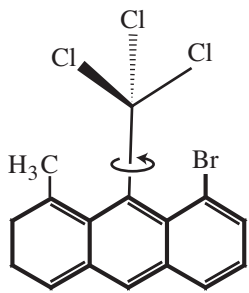}
\caption{The point group of this molecule is the cyclic group $\mathbb{Z}_2$.}
\label{Z3}
\end{figure}

Knowing the symmetries of a molecule helps to predict its chemical behavior.
But, what exactly do we mean by ``symmetries?''  If we consider only rigid molecules, then the molecular symmetries are
rotations, reflections, and reflections composed with rotations. Chemists have defined
the {\it point group} of a molecule as its group of rigid symmetries.  While this is a useful tool for rigid molecules, it can be misleading when applied to non-rigid molecules.  For example, consider the molecule illustrated in Figure \ref{Z3}.  The Cl on the far left is in front of the page, the Cl in the middle is behind the page, and the Cl on the right is in the page.

The point group of the molecule in Figure 1 is the cyclic group $\mathbb{Z}_2$ because the only rigid symmetry of the molecule is a reflection through the plane which contains the three hexagons.  However, the three Cl's at the top can rotate around the bond connecting them to the hexagons (as indicated by the arrow) while the hexagons themselves remain fixed. We would like to define a group which includes both the rigid and non-rigid symmetries of molecules.  

We begin by defining a {\it molecular graph} as a graph embedded in $\mathbb{R}^3$ where the edges represent bonds and the vertices represent atoms or groups of atoms.  We then use the following definition.

\begin{definition}
An \emph{automorphism} of a molecular graph $\gamma$ is a permutation of the vertices preserving adjacency, taking atoms of a given type to atoms of the same type (e.g., carbons go to carbons, oxygens go to oxygens, and so on).  $\mathrm{Aut}(\gamma)$ is defined as the group of automorphisms of $\gamma$.
\end{definition}

The automorphism group is the group of symmetries of the abstract graph, independent of any embedding in $\mathbb{R}^3$.  However the symmetries of a molecule are only those symmetries of the abstract graph which occur as symmetries of the molecular graph.   In particular, we have the following definition.

\begin{definition}We define the \emph{molecular symmetry group} of a molecular graph $\gamma$ as the subgroup of $\mathrm{Aut}(\gamma)$ induced by chemically possible motions taking the molecule to itself together with reflections which take the molecule to itself.  
\end{definition}

The molecular symmetry group of the graph illustrated in Figure \ref{Z3}  is $D_3$ (the dihedral group with six elements) and is generated by the automorphism induced by a planar reflection together with the automorphism induced by a $120^{\circ}$ rotation of the three Cl's at the top keeping the rest of the graph fixed.   

While the molecular symmetry group makes sense chemically, it cannot be defined mathematically because the existence of a particular molecular symmetry may depend on the flexibility of the molecule or the possibility of rotating around certain bonds within the molecule.  Rather than assuming all molecules are completely rigid as the point group does, we now consider the group obtained by treating all molecules as if they were completely flexible.

\begin{definition}
The \emph{topological symmetry group} of a graph $\Gamma$ embedded in $\mathbb{R}^3$ is the subgroup of $\Aut(\Gamma)$ induced by diffeomorphisms of the pair $(\mathbb{R}^3, \Gamma)$.  It is denoted by $\mathrm{TSG}(\Gamma)$
\end{definition}

\begin{figure} [htp]
\includegraphics{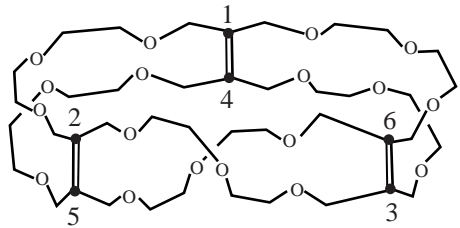}
\caption{The topological symmetry group of this molecule is $D_6$.}
\label{mobius}
\end{figure}

Consider the graph of the molecular M\"{o}bius ladder illustrated in Figure \ref{mobius}.  This molecule is large enough to be somewhat flexible. We number some of the vertices so that we can write the automorphisms more conveniently.  The automorphism $(23)(56)(14)$ is induced by turning the molecule upside down.  This is the only non-trivial automorphism which is induced by a rigid symmetry.  Thus the point group of the molecule is $\mathbb{Z}_2$.   However, the automorphism $(123456)$ is induced by rotating the molecule by $120^{\circ}$ while slithering the half-twist back to its original position.
Thus the topological symmetry group of this molecule is $D_6$ (the dihedral group of order 12).  Because of the flexibility of the molecule, this is also the molecular symmetry group.

While the motivation for the study of topological symmetry groups came from considering the symmetries of non-rigid molecules, we can consider the topological symmetry group of any graph embedded in $\mathbb{R}^3$.  In fact, the study of topological symmetry groups is a natural extension of the study of the symmetries of knots.  Thus, as is typical in knot theory, we prefer to consider graphs embedded in $S^3$ rather than in $\mathbb{R}^3$.  Observe that the topological symmetry group of a graph embedded in $\mathbb{R}^3$ is the same as it is if we consider the same embedding of the graph in $S^3$.

\begin {definition}
Let $\Gamma$ be a graph embedded in $S^3$. The \emph{ orientation preserving topological symmetry group}, $\mathrm{TSG}_+(\Gamma)$, is the subgroup of $\mathrm{TSG}(\Gamma)$ induced by orientation preserving diffeomorphisms of the pair $(S^3, \Gamma)$.  \end{definition}

 Observe that either $\mathrm{TSG}_+(\Gamma)=\mathrm{TSG}(\Gamma)$ or $\mathrm{TSG}_+(\Gamma)$ is a normal subgroup of $\mathrm{TSG}(\Gamma)$ of index 2.  So understanding $\TSG(\Gamma)$ is a key step in understanding $\mathrm{TSG}(\Gamma)$.  For the rest of this paper we will focus our attention on the orientation preserving topological symmetry group.  However, for the sake of simplicity we will abuse terminology and refer to these groups simply as topological symmetry groups.

The general question we are interested in is: for a given graph, what groups can occur as topological symmetry groups?   In \cite{FNPT} it was shown that not every finite group can occur as the topological symmetry group of some graph in $S^3$.  In particular, $\mathrm{TSG}_+(\Gamma)$ cannot be the alternating group $A_n$ for $n>5$ for any embedded graph $\Gamma$ in $S^3$.  However, there is no known classification of all possible topological symmetry groups of graphs in $S^3$.  A {\it complete graph}, $K_n$, is a graph with $n$ vertices and an edge between every pair of vertices.  The class of complete graphs is an interesting family of graphs to consider because $\Aut(K_n)$ is the symmetric group $S_n$, which is the largest automorphism group of any graph with $n$ vertices. In \cite{FNT}, Flapan, Naimi, and Tamvakis proved the following theorem, characterizing which finite groups can occur as topological symmetry groups of embeddings of complete graphs in $S^3$.

\medskip \begin{complete} \cite{FNT}
\label{TSG2} 
A finite group $H$ is $\TSG(\Gamma )$ for some embedding $\Gamma$ of 
a complete graph in $S^3$ if and only if $H$ is 
isomorphic to a finite cyclic group, a dihedral group, $S_4$, $A_4$, $A_5$, or
a subgroup of $D_{m}\times D_{m}$ for some odd 
$m$.
\end{complete}

 Observe that the Complete Graph Theorem does not tell us for a given $n$, what groups can occur as the topological symmetry group of some embedding of $K_n$ in $S^3$.  In this paper we characterize what groups can occur as the topological symmetry group of an embedding of any complete graph of the form $K_{4r+3}$ in $S^3$.

 \section{Topological symmetry groups of complete graphs}

In addition to the Complete Graph Theorem, we will make use of several results from other papers.  The following result from  \cite{FNPT} shows us that for 3-connected graphs (i.e., those which cannot be disconnected or reduced to a single vertex by deleting fewer than 3 vertices), we only need to consider topological symmetry groups that are induced by a finite subgroup of $\Diff(S^3)$ (i.e., the group of orientation preserving diffeomorphisms of $S^3$).  
  
  \begin{finiteness}\cite{FNPT}
  \label{TSG1}
  Let 
$H=\TSG(\Gamma )$ for some 
$3$-connected graph $\Gamma $ embedded in $S^3$.  Then $\Gamma $ 
can be re-embedded  as $\Delta $ such that $H\leq\TSG(\Delta )$ and 
$\TSG(\Delta )$ is induced by an isomorphic finite
subgroup of $\Diff(S^{3})$.
\end{finiteness}

We contrast this result with the example, illustrated in Figure \ref{Z3}, in which the rotation of the three Cl's cannot be induced by any finite order diffeomorphism.   Observe that the graph in Figure 1 can be disconnected by removing a single vertex, and hence it is not 3-connected.

The following theorem classifies which automorphisms of $K_n$ can be induced by a finite order orientation preserving diffeomorphism of some embedding of $K_n$ in $S^3$.  For simplicity we will use the word {\it cycle} only for non-trivial cycles of an automorphism.  Note that, if $n=4r+3$, then $n>6$.

\begin{auto} \cite{Fl}
\label{automorphism}
Let 
$K_{n}$ be a complete graph on $n>6$ vertices 
and let $\varphi $ be an automorphism of $K_{n}$.  Then there 
is an embedding $\Gamma$ of $K_{n}$ in $S^{3}$ such 
that $\varphi $ is induced by an orientation 
preserving diffeomorphism $h$ of $(S^3,\Gamma)$ of order $m$ if 
and only if the cycles and fixed vertices of $\varphi$ can be described by one of the following: 
\begin{enumerate}
\item $m>2$ is even, all 
cycles of $\varphi $ are of order $m$, and $\varphi$ fixes no vertices.

\item  $m=2$, all
cycles of $\varphi $ are of order $m$, and $\varphi$ fixes 
at most two vertices.

\item  $m$ is odd, all 
 cycles of $\varphi $ are of order $m$, 
and $\varphi$ fixes at most three vertices.

\item $m$ is an odd multiple of $3$ 
and $m>3$, all cycles of $\varphi $ are of order 
$m$ except one of order $3$, and $\varphi$ fixes no vertices.
\end{enumerate}
\end{auto}

We will say that an automorphism $\varphi$ of $K_n$ is of {\it type 4}, if it is described by Condition (4) of the Automorphism Theorem.  

We now prove some general lemmas about automorphisms of graphs and  subgroups of $\Diff(S^3)$ which are the product of two cyclic groups.  We will use $\fix(h)$ to denote the fixed point set of a diffeomorphism $h$ of $S^3$.  Note that if $h$ is orientation preserving, then by Smith Theory \cite{Sm} either $\fix(h)\cong S^1$ or $\fix(f)=\emptyset$.

\bigskip

\begin{orbits} \label{orbit lemma}
Suppose $\alpha$ and $\beta$ are commuting automorphisms of a finite set $V_1$.  Then $\beta$ takes $\alpha$-orbits to $\alpha$-orbits of the same length.
\end{orbits}

\begin{proof}

Let $A$ denote a minimal length $\alpha$-orbit (i.e., no $\alpha$-orbit has a shorter length). Then

\[ \beta(A) = \beta(\alpha(A)) = \alpha(\beta(A)) \].

\noindent Thus $\alpha$ takes $\beta(A)$ to $\beta(A)$.  This implies that $\beta(A)$ is a union of $\alpha$-orbits.  However, $|\beta(A)| = |A|$, which is a minimal length $\alpha$-orbit.  Hence $\beta(A)$ is actually a single $\alpha$-orbit.

Let $A_1, A_2, ..., A_n$ denote all the minimal length $\alpha$-orbits in $V_1$, then for each $i$,  there exists a $j$ such that $\beta(A_i) =A_j$.  Thus
\[ \beta(A_1 \cup A_2 \cup ... \cup A_n) = A_1 \cup A_2 \cup ... \cup A_n \]

Let $V_2 = V_1 - (A_1 \cup A_2 \cup ... \cup A_n)$, then $\alpha(V_2) = V_2$ and $\beta(V_2) = V_2$.  Now start with $V_2$, and repeat the above argument as necessary to see that  $\beta$ takes every $\alpha$-orbit to an $\alpha$-orbit of the same length.  \end{proof}

\bigskip

\begin{disjoint}\label{disjoint}
Suppose $g, h \in \mathrm{Diff}_+(S^3)$ such that $\langle g, h\rangle = \mathbb{Z}_p \times \mathbb{Z}_q$ is not cyclic or equal to $D_2$.  Then $\mathrm{fix}(g)$ and $\mathrm{fix}(h)$ are disjoint.
\end{disjoint}

\begin{proof}
Let $G = \langle g, h \rangle$.  Suppose there exists $x \in \mathrm{fix}(g) \cap \mathrm{fix}(h)$.  Let $B(x)$ denote a regular neighborhood of $x$ in $S^3$, and let $N(x)=\bigcap_{f\in G}f(B(x))$.  Since $G$ is finite and fixes $x$, $N(x)$ is a ball around $x$ which is setwise invariant under $G$.

Let $H$ be the restriction of $G$ to the sphere, $\partial N(x)$.  By Smith Theory \cite{Sm}, if some $f \in G$ pointwise fixes $\partial N(x)$, then $f$ must be the identity.  Thus $H \cong G=\mathbb{Z}_p \times \mathbb{Z}_q$, and $H$ is neither cyclic nor equal to $D_2$.  Since this is impossible, $\mathrm{fix}(g)$ and $\mathrm{fix}(h)$ must in fact be disjoint.  \end{proof}

\bigskip

In the next several lemmas we put restrictions on the types of automorphisms that can be induced on an embedding $\Gamma$ of $K_n$ by a group $\mathbb{Z}_p \times \mathbb{Z}_q\leq\Diff(S^3)$ 

\bigskip

\begin{3cycle}
Let $\Gamma$ be an embedding of $K_n$ in $S^3$, and let $\alpha$ and $\beta$ be orientation preserving diffeomorphisms of $(S^3, \Gamma)$ of odd orders p and q respectively such that $p$, $q>1$ and $p|q$.  Suppose that $\langle \alpha, \beta \rangle = \mathbb{Z}_p \times \mathbb{Z}_q$, and $\langle \alpha \rangle \cap \langle \beta \rangle = \langle e \rangle$.  Then the following hold.
\begin{enumerate}

\item $\alpha$ cannot be of type 4, and if $\beta$ is of type 4 then $p=3$.  

\item If $\alpha$ and $\beta$ have no fixed vertices and either $\alpha$ or $\beta$ has no 3-cycles, then $pq|n$.

\item There are at most two disjoint sets of 3 vertices which are setwise invariant under both $\alpha$ and $\beta$.

\end{enumerate}

 \end{3cycle}
\medskip

\begin{proof} 

Suppose that there is some diffeomorphism $\delta \in \langle \alpha, \beta \rangle$ which is of type 4.  Then, by the Automorphism Theorem, there exists a unique $\delta$-cycle $A$ of order 3.  Now by the Orbits Lemma, $\alpha(A) = A$ and $\beta(A) = A$.  Thus $\alpha^3$ and $\beta^3$ fix $A$ pointwise, and therefore by the Disjoint Fixed Points Lemma, at least one of $p=3$ or $q=3$.  If $\alpha$ is of type 4, then $p\not=3$ and since $p|q$, $q\not =3$.  Hence $\alpha$ cannot be of type 4.  If $\beta$ is of type 4, then $q\not=3$ and hence $p=3$.   Thus Conclusion (1) holds.

Suppose $\alpha$ and $\beta$ have no fixed vertices.  Let $A$ and $B$ denote the $\alpha$-orbit and $\beta$-orbit respectively of some vertex $v$.  Suppose that for some $i<q$, $\beta^i(A)=A$.  Then $\beta^i(v) = \alpha^k(v)$ for some $k < p$.   Let $w\in A$. Then $w = \alpha^s(v)$ for some $s < p$.  Thus 
$$\alpha^k(w) = \alpha^k(\alpha^s(v)) = \alpha^s(\alpha^k(v)) = \alpha^s(\beta^i(v)) = \beta^i(\alpha^s(v)) =\beta^i(w)$$

Hence $\alpha^k \beta^{-i}$ fixes every $w\in A$.  Also, since $\langle \alpha \rangle \cap \langle \beta \rangle = \langle e \rangle$, $\alpha^k\beta^{-i}$ is not the identity. By the Automorphism Theorem, this implies that $|A|\leq3$.  Since $p$ is odd and $v$ is not fixed by $\alpha$, we must have $|A|=3$.  Similarly, if for some $j<p$, $\alpha^j(B)=B$, then $|B|=3$.

 Suppose that $\alpha$ has no 3-cycles. Then by the above argument, for every $i<q$, $\beta^i(A)\not =A$.   We know by the Orbits Lemma that $\beta$ takes $\alpha$-orbits to $\alpha$-orbits. Since $\beta$ has order $q$, $\beta$ permutes the $\alpha$-orbits in cycles of length $q$.  Since $\alpha$ has no fixed vertices or 3-cycles, every $\alpha$-orbit has length $p$.  This implies that $pq|n$.  By switching the roles of $\alpha$ and  $\beta$, we obtain the same result when $\beta$ has no 3-cycles. Thus Conclusion (2) follows.

Now suppose that $A_1$, $A_2$, and $A_3$ are disjoint sets of three vertices such that for each $i$, $\alpha(A_i)=A_i$ and $\beta(A_i)=A_i$.
Then $\alpha$ and $\beta$ each have at least two 3-cycles and either a third 3-cycle or 3 fixed vertices.  Hence by the Automorphism Theorem $p=3=q$. For each $i$, we take the union of the three edges of $\Gamma$ joining pairs of vertices in $A_i$ to obtain a triangle $B_i$ which is setwise invariant under both $\alpha$ and $\beta$.    Thus $B_1$, $B_2$, and $B_3$ are disjoint simple closed curves which are each setwise invariant under the group $\langle\alpha,\beta\rangle\leq \Diff(S^3)$.

Observe that no pair of $\alpha$, $\beta$, $\alpha \beta$, and $\alpha \beta^2$ are equal or inverses of one another.  Since all non-trivial elements of $\Z_3 \times \Z_3$ are of order 3, for each $i$, there are at most three distinct diffeomorphisms leaving $B_i$ setwise invaraint.  Thus for each $i$, at least two of $\alpha$, $\beta$, $\alpha \beta$, and $\alpha \beta^2$ induce the same diffeomorphism of $B_i$.  Hence for each $i$, there is a non-trivial $f_i\in \langle\alpha,\beta\rangle$ which pointwise fixes $B_i$.  

By the Geometrization Conjecture for Orbifolds \cite{BLP}, the group of diffeomorphisms $\langle\alpha,\beta\rangle$ is conjugate to a group of isometries of $S^3$.  However, up to conjugacy there is only one group $H$ of isometries of $S^3$ which is isomorphic to $\mathbb{Z}_3 \times \mathbb{Z}_3$.  The group $H$ consists of all isometries of the unit sphere $S^3$ of $\mathbb{C}^2$ of the form $(z_1,z_2)\mapsto(\omega^{k_1}z_1,\omega^{k_2}z_2)$ where $\omega$ is a third root of unity.  The only possible simple closed curves which can be the fixed point set of a non-trivial element of this group of isometries, are the intersection of $S^3$ with one of the two axes in $\mathbb{C}^2$.  Therefore at most two simple closed curves are setwise invariant under $\langle\alpha,\beta\rangle$.  Thus there cannot be non-trivial elements of $\langle\alpha,\beta\rangle$ which pointwise fix each $B_i$.  It follows that there are at most two sets of 3 vertices which are setwise invariant under both $\alpha$ and $\beta$, and hence Conclusion (3) follows. \end{proof}

\bigskip

\begin{fixed} 
Let $\Gamma$ be an embedding of $K_n$ in $S^3$.  Let $\alpha$ and $\beta$ be orientation preserving diffeomorphisms of $(S^3, \Gamma)$ of odd orders p and q respectively such that $p$, $q>1$ and $p|q$.   Suppose that $\langle \alpha, \beta \rangle = \mathbb{Z}_p \times \mathbb{Z}_q$ such that $\langle \alpha \rangle \cap \langle \beta \rangle = \langle e \rangle$.  Then the following are true.
\begin{enumerate}
\item If either $\alpha$ or $\beta$ fixes any vertices, then it fixes 3 vertices. 

\item If $\beta$ fixes any vertices then $p=3$, and if $\alpha$ fixes any vertices, then $p=q=3$. 
\end{enumerate}
\end{fixed}

\begin{proof}
Suppose that $\alpha$ fixes precisely one vertex $v_1$.  Then by the Orbits Lemma, $\beta(v_1)=v_1$, since $\{v_1\}$ is the only $\alpha$-orbit of length 1.  This implies $v_1$ is a fixed vertex of $\beta$.  However, since $\langle \alpha, \beta \rangle$ is neither cyclic nor $D_2$, we can apply the Disjoint Fixed Points Lemma to get a contradiction.  Thus neither $\alpha$ nor $\beta$ can have precisely one fixed vertex.  Suppose that $\alpha$ fixes precisely two vertices, $v_1$ and $v_2$. Then by the Orbits Lemma, $\beta(\{v_1,v_2\})=\{v_1,v_2\}$,  and by the Disjoint Fixed Points Lemma, $\beta$ cannot fix either $v_i$.  This implies that $\beta$ has an orbit of length 2, which contradicts the Automorphism Theorem since $p$ and $q$ are of odd order.  Since we did not use the hypothesis that $p|q$ (except in assuming that $\langle \alpha, \beta \rangle$ is not cyclic), the roles of $\alpha$ and $\beta$ could be switched.  Thus Conclusion (1) holds.

Suppose that $\beta$ has 3 fixed vertices. Then by the Orbits and Disjoint Fixed Points Lemmas, these 3 vertices are a 3-cycle of $\alpha$. Now $\alpha$ cannot be of type 4 by the 3-Cycle Lemma.  Thus $p=3$.  Suppose that $\alpha$ has 3 fixed vertices.  Then these vertices form a 3-cycle of $\beta$.  Now $\beta^3$ and $\alpha$ have these 3 fixed vertices in common. Thus  by the Disjoint Fixed Points lemma $\beta^3$ must be the identity, and hence $q=3$ which implies that $p=3$.  Therefore, Conclusion (2) holds.  \end{proof}
\bigskip

\begin{pq} \label{pq}
Let $\Gamma$ be an embedding of $K_n$.  Let $\alpha$ and $\beta$ be orientation preserving diffeomorphisms of $(S^3, \Gamma)$ of odd orders p and q respectively such that $p$, $q>1$ and $p|q$.  Suppose that $\langle \alpha, \beta \rangle = \mathbb{Z}_p \times \mathbb{Z}_q$ such that $\langle \alpha \rangle \cap \langle \beta \rangle = \langle e \rangle$.  Then the following conclusions hold:

\begin{enumerate}
\item If $p > 3$, then $pq|n$.
\item If $p=3$ and $\beta$ is of type 4, then $pq|n-3$.
\item If $p=3$, $q \neq 3$, and $\beta$ is not of type 4, then either $pq|n$ or $pq|n-3$.
\item If $p=q=3$, then either $pq|n$ or $pq|n-3$ or $pq|n-6$.
\end{enumerate}
\end{pq}

\begin{proof}  Suppose that $p>3$.  Then by the Fixed Vertex Lemma, neither $\alpha$ nor $\beta$ fixes any vertices.  By the 3-Cycle Lemma, $\alpha$ cannot be of type 4.  Hence $\alpha$ has no 3-cycles, and hence by the 3-Cycle Lemma $pq|n$.  Thus we obtain Conclusion (1).  For the rest of the proof, we assume that $p=3$.

Suppose $\beta$ is of type 4.  By the Automorphism Theorem, $\beta$ has precisely one 3-cycle A, and $\beta$ fixes no vertices.  Also, since $q \neq 3$, by the Fixed Vertex Lemma $\alpha$ fixes no vertices.  By the Orbits Lemma, $\alpha$ takes $\beta$-orbits to $\beta$-orbits of the same length.  Hence, $\alpha(A) = A$.  Let $\Gamma'$ denote the embedding of $K_{n-3}$ obtained from $\Gamma$ by deleting the vertices in $A$ and the edges containing them.  Then $\alpha$ and $\beta$ leave $\Gamma'$ setwise invariant, neither $\alpha$ nor $\beta$ fixes any vertices of $\Gamma'$, and $\beta$ has no 3-cycles in $\Gamma'$. Thus by the 3-Cycle Lemma applied to $\Gamma'$, we know that $pq|n-3$.   Hence Conclusion (2) holds.

Next suppose that $q \neq 3$ and $\beta$ is not of type 4.  Then $\beta$ has no orbits of length 3.  By the Fixed Vertex Lemma, $\alpha$ has no fixed vertices and $\beta$ has either 0 or 3 fixed vertices.  If $\beta$ has no fixed vertices, then by the 3-Cycle Lemma, $pq|n$.  Suppose that fix($\beta$)=$A$ contains 3 vertices.  By the Orbits lemma $\alpha(A) = A$. Define $\Gamma'$ to be the embedding of $K_{n-3}$ obtained from $\Gamma$ by removing the vertices in $A$ and the edges containing them.  Then $\alpha$ and $\beta$ leave $\Gamma'$ setwise invariant, neither $\alpha$ nor $\beta$ fixes any vertices of $\Gamma'$, and $\beta$ has no 3-cycles.   Thus we can apply the 3-Cycle Lemma to $\Gamma'$ to conclude that $pq |n-3$.  Thus Conclusion (3) holds.

Finally, suppose $p=q=3$. By the Fixed Vertex Lemma,  $\alpha$ and $\beta$ each have either 0 or 3 fixed vertices.  By the Disjoint Fixed Points Lemma, the sets of fixed vertices of $\alpha$ and $\beta$ are disjoint, and by the Orbits Lemma the fixed vertices of one of $\alpha$ and $\beta$ are setwise invariant under the other.   Let $X$ denote the union of all of the sets of 3 vertices which are setwise invariant under both $\alpha$ and $\beta$.  Thus $X$ includes the fixed vertices of $\alpha$ and $\beta$, as well as any set of 3 vertices which is an orbit of both $\alpha$ and $\beta$.  By the 3-Cycle Lemma, no more than two disjoint sets of 3 vertices are setwise invariant under both $\alpha$ and $\beta$.  Thus $|X|=$ 0, 3, or 6.

Let $\Gamma''$ be obtained from $\Gamma$ by removing the vertices in $X$ together with the edges between them.  Thus $\alpha$ and $\beta$ leave $\Gamma''$ setwise invariant and no set of 3 vertices of $\Gamma''$ are setwise invariant under both $\alpha$ and $\beta$.  Furthermore, no vertices of $\Gamma''$ are fixed by $\alpha$, and hence all $\alpha$-orbits have $p=3$ vertices.   We know by the Orbits Lemma that $\beta$ takes $\alpha$-orbits to $\alpha$-orbits; and no orbits of $\alpha$ are also orbits of $\beta$.  Therefore, $\beta$ permutes the $\alpha$-orbits in cycles of length $q=3$.   Let $m$ denote the number of vertices of $\Gamma''$, then $pq|m$.   However, $m$ is either $n$, $n-3$, or $n-6$.  Thus Conclusion (4) follows.\end{proof}

\section{Topological symmetry groups of $K_{4r+3}$}

Now we focus on embeddings of complete graphs of the form $K_{4r+3}$.

\begin{D2} 
There is no embedding $\Gamma$ of any $K_{4r+3}$ in $S^3$ such that $D_2 \leq \mathrm{TSG}_+(\Gamma)$.
\end{D2}

\begin{proof}
 Suppose that $\Gamma$ is an embedding of some $K_{4r+3}$ such that $D_2 \leq \mathrm{TSG}_+(\Gamma)$. Let $\varphi_1$ and $\varphi_2$ be distinct non-trivial elements of $D_2 \leq \mathrm{TSG}_+(\Gamma)$.  Since the number of vertices of $\Gamma$ is $4r +3$, by the Automorphism Theorem, each $\varphi_i$ must be composed of $(2r +1)$ 2-cycles with precisely one fixed vertex.  Thus, $\varphi_1\varphi_2$ can be written as the product of $(4r +2)$ (not necessarily disjoint) 2-cycles.  However, $\varphi_1\varphi_2$ also has order 2, and hence by the Automorphism Theorem can also be written as a product of $(2r+1)$ 2-cycles.  However, no automorphism can be written as both the product of an even number of 2-cycles and the product of an odd number of 2-cycles.  Hence such an embedding of $K_{4r+3}$ cannot exist.\end{proof}

\medskip

Observe that $D_2 $ is contained in the groups $A_4$, $S_4$, and $A_5$.  Thus by the $D_2$ lemma, there is no embedding $\Gamma$ of $K_{4r+3}$ in $S^3$ such that $\mathrm{TSG}_+(\Gamma)$ is $A_4$, $S_4$, or $A_5$.  Now it follows from the Complete Graph Theorem, that any embedding $\Gamma$ of $K_{4r+3}$ in $S^3$, $\mathrm{TSG}_+(\Gamma)$ is either cyclic, dihedral, or a subgroup of $D_m \times D_m$ for some odd $m$.  The following is our main result.

\bigskip

\begin{4r+3}  Let $n=4r+3$.  A finite group $G$ is isomorphic to $\TSG(\Gamma )$ for some embedding $\Gamma$ of 
$K_n$ in $S^3$ if and only if 
one of the following holds:
\begin{enumerate}

\item $G$ is $\mathbb{Z}_2$
\item $G$ is $D_p$ or $\mathbb{Z}_p$ where $p$ is odd and either $p|n$, $p|n-1$, $p|n-2$, or $p|n-3$.
\item $G=\mathbb{Z}_p \times \mathbb{Z}_q$ where $p$ and $q$ are odd, $p|q$ and one of the following holds.
\begin{itemize}

\item $pq|n$.

\item $p=3$ and $pq|n-3$.
\item $p=q=3$ and $pq|n-6$.
\end{itemize}
\end{enumerate}
\end{4r+3}
\bigskip

In order to construct the required embeddings, we will first define an embedding of the vertices of $K_{4r+3}$ which is setwise invariant under a finite subgroup of $\Diff(S^3)$ and then use the following result to embed the edges.

\begin{edge}  \cite{FMN2} Let $G$ be a finite subgroup of $\Diff(S^3)$ which has the property that for any even order $g\in G$, $\fix(g)\not =\emptyset$ if and only if  $\mathrm{order}(g)=2$.  Let $\gamma$ be a graph whose vertices are embedded in $S^3$ as a set $V$ which is invariant under $G$ such that $G$ induces a faithful action of $\gamma$. Suppose that the vertices of $\gamma$ satisfy the following additional hypotheses:
\begin{enumerate}

\item No pair of adjacent vertices is equal to $\fix(h)\cap \fix(g)$ for some $h$, $g\in G$.

\item At most one pair of adjacent vertices is fixed by an order 2 element of $G$, and no pair of adjacent vertices is interchanged by any $g\in G$ with $\mathrm{order}(g)\not =2$.

\item Any pair of adjacent vertices $\{v, w\}$ which is fixed by a non-trivial $g\in G$ bounds an arc $A_{vw}\subseteq \fix(g)$ whose interior is disjoint from $V$ and if the set $\{v,w\}$ or a point of $\INT(A_{vw})$ is invariant under some $f\in G$ then $f(A_{vw})=A_{vw}$.
\end{enumerate}

Then there is an embedding of the edges of $\gamma$ in $S^3$ such that the resulting embedding of $\gamma$ is setwise invariant under $G$.  
\end{edge}

The Edge Embedding Lemma will give us an embedding $\Gamma$ of $K_{4r+3}$ in $S^3$ with $G\leq \TSG(\Gamma)$.  In order to create an embedding $\Gamma'$ with $G= \TSG(\Gamma')$ we will need the following result which was proved as Theorem 2 of \cite{FMN1}.

\begin{subgroup}\cite{FMN1}
Let $n > 6$, and suppose that $\Gamma$ is an embedding of $K_n$ in $S^3$ such that $\TSG(\Gamma)$ is cyclic, dihedral, or a subgroup of $D_m \times D_m$ for some odd $m$.  Then for every $H\leq \TSG(\Gamma)$, there is an embedding $\Gamma'$ of $K_n$ such that $H = \TSG(\Gamma')$.
\end{subgroup}
\medskip

With these results in hand, we prove our theorem as follows.
\medskip
   
\begin{proof}  We begin by assuming that $G=\TSG(\Gamma)$ for some embedding $\Gamma$ of $K_n$ in $S^3$.   As we observed above, it follows from the Complete Graph Theorem and the $D_2$ Lemma that $G$ is either a finite cyclic group, a dihedral group, or a subgroup of $D_m \times D_m$ for some odd $m$.  
By the Finiteness Theorem, there exists a re-embedding $\Delta$ of $K_n$ in $S^3$ such that $G$ is induced by an isomorphic finite subgroup of $\Diff(S^3)$.

Suppose $G$ is equal to $\mathbb{Z}_p$ or $D_p$ for some integer $p$. If $p=2$, then $G=\mathbb{Z}_2$ since we know by the $D_2$ Lemma that $G\not=D_2$.  Thus we assume that $p>2$.  Now suppose that $p$ is even.  Since  $n=4r+3$, $n>6$.  
Thus by the Automorphism Theorem, the elements of $G$ of order $p$ are composed of $p$-cycles with no fixed vertices.  Hence $p|n$.  However, since $p$ is even and $n=4r+3$, this is impossible.  Thus $p$ must be odd.  Now by the Automorphism Theorem, the elements of $G$ are each composed entirely of $p$-cycles with up to 3 fixed vertices.  This means that either $p|n$, $p|n-1$, $p|n-2$, or $p|n-3$.

Next we assume that $G \leq D_m \times D_m$ for some odd $m$, but $G$ is neither cyclic nor dihedral and $G$ does not contain $D_2$.  It follows that there exist $\alpha$, $\beta\in G$ of odd order $a$ and $b$ respectively such that $\langle\alpha,\beta\rangle=\mathbb{Z}_a\times\mathbb{Z}_b\leq G$ and $\mathbb{Z}_a\times\mathbb{Z}_b$ is not cyclic.  

Suppose for the sake of contradiction that $G$ contains an order 2 element $\varphi$.   Since $G \leq D_m \times D_m$, for every $g\in G$, either $g\varphi=\varphi g$ or $g\varphi=\varphi g^{-1}$.  If $\alpha \varphi=\varphi\alpha$, then $\alpha\varphi$ has order $2a$, since $a$ is odd.  However, since $2a>2$, by the Automorphism Theorem, $\alpha\varphi$ fixes no vertices.  But this is impossible since $2a$ does not divide $4r+3$.  Thus we must have $\alpha \varphi=\varphi\alpha^{-1}$.   Now, since $\varphi$ has order 2 and $n=4r+3$, $\varphi$ fixes some vertex $v$ of $\Gamma$.  Thus $\varphi\alpha^{-1}(v)=\alpha\varphi(v)=\alpha(v)$, and hence $(\varphi\alpha^{-1})^2(v)=\alpha^{2}(v)$.  However, since $(\varphi\alpha^{-1})^2(v)=\varphi\alpha^{-1}\varphi\alpha^{-1}(v)=\varphi^2\alpha\alpha^{-1}(v)=v$.  Thus $\alpha^2$ is a non-trivial element of $\mathbb{Z}_a\times\mathbb{Z}_b$ which fixes $v$. Using an analogous argument we see that $\beta^2$ is a non-trivial element of $\mathbb{Z}_a\times\mathbb{Z}_b$ which also fixes $v$.  However, since $\alpha$ and $\beta$ both have odd order, $\langle\alpha^2,\beta^2\rangle=\langle\alpha,\beta\rangle=\mathbb{Z}_a\times\mathbb{Z}_b$, which is not cyclic or equal to $D_2$.  Thus we can apply the Disjoint Fixed Points Lemma to get a contradiction.  Thus $G$ does not contain any element of order 2.  

It follows that $G=\mathbb{Z}_p\times\mathbb{Z}_q$ for some odd $p$ and $q$ such that $p|q$.  Now the required conclusions all follow from the pq Lemma.
\medskip

In order to prove the converse, we consider two cases as follows.  
\medskip

\noindent {\bf Case 1:}  Either $G=D_p$ or $\Z_p$ for $p$ odd, or $G=\Z_2$

First let $G=D_p$ with $p$ odd and let $G'\cong D_p$ be a subgroup of $\Diff(S^3)$ generated by a rotation $g$ by $\frac{2\pi}{p}$ about a circle  $C_g$ together with a rotation $f$ by $\pi$ about a circle $C_f$ which meets $C_g$ in two points.  Observe that the orbit of $C_f$ under $G'$ is $p$ circles which are each the fixed point set of an order 2 element of $G'$.

  Suppose that $n=kp$ for some integer $k$.  First observe that since $n=4r+3$ and $p$ are both odd, $k$ must be odd.    Let $x_1$ be a point on one of the arcs in $C_f-C_g$.  The orbit of $x_1$ consists of $p$ points $x_1$, $x_2$, \dots, $x_p$, each on the fixed point set of a distinct element of order 2 in $G'$.  We embed $p$ of the $n$ vertices as the points $x_1$, $x_2$, \dots, $x_p$.  
  
   Let $Y$ denote the union of the fixed point sets of all non-trivial elements of $G'$.  Let $B$ denote a ball which is disjoint from $Y$ and from its image under every non-trivial element of $G'$.  Since $k$ is odd, we can embed $\frac{k-1}{2}$ points in $B$.  Since the order of $G'$ is $2p$, the orbit of these $\frac{k-1}{2}$ points is $(k-1)p$ points $z_1$, \dots $z_{p(k-1)}$ in $S^3-Y$.  We embed the remaining $(k-1)p$ vertices as this set of points.  Observe that the sets of embedded vertices $X=\{x_1,\dots, x_p\}$ and $Z=\{z_1,\dots, z_{p(k-1)}\}$ are setwise invariant under $G'$, and $G'|X\cup Z$ induces $G$ on $K_n$.

Since no pair of vertices is fixed by a non-trivial element of $G'$ and $G'$ contains no even order elements with order greater than 2, it is not hard to check that the hypotheses of the Edge Embedding Lemma are satisfied for the embedded vertices of $K_n$.  Thus we can apply the Edge Embedding Lemma to get an embedding $\Gamma$ of $K_n$ with $D_p\leq \TSG(\Gamma)$.  Finally, since $n=4r+3$, we know by the forward direction of our proof that $\TSG(\Gamma)$ is either cyclic, dihedral, or a subgroup of $D_m \times D_m$ for some odd $m$.  Thus we can apply the Subgroup Theorem to get the required embedding $\Gamma'$ of $K_n$ with $D_p= \TSG(\Gamma')$.

Now let $n=kp+2$ for some $k$.  Since $n=4r+3$ and $p$ are both odd, $k$ is again odd. Thus we can embed $kp$ vertices as described above.  Now we add two additional vertices $v$ and $w$  in $C_g-C_f$ so that $f$ interchanges them.   Thus we have embedded all of the vertices of $K_{kp+2}$ as a set $V$ such that $G'|V$ induces $G$ on $K_{kp+2}$.  Observe that $v$ and $w$ are the only pair of vertices which are fixed by a non-trivial element of $G$.  Let $A_{vw}$ be one of the arcs in $C_g$ which is bounded by $v$ and $w$.  Then $A_{vw}$ satisfies hypothesis (3) of the Edge Embedding Lemma.  Thus, as above, we apply the Edge Embedding Lemma and then the Subgroup Lemma to get the required embedding $\Gamma'$ of $K_n$ with $D_p= \TSG(\Gamma')$.  

Next let $n=kp+1$ for some integer $k$.  In this case, since $n=4r+3$ and $p$ are both odd, $k$ must be even.  Let $B$ be the ball described above, and choose $\frac{k}{2}$ points in $B$. The orbit of these points under $G'$ consists of $kp$ points which will be embedded vertices.  We embed the final vertex $x$ as one of the points in $C_g\cap C_f$.  Thus $x$ is fixed by every element of $G'$. Since again no pair of vertices is fixed by a non-trivial element of $G'$, we again obtain the required embedding of $K_n$ by applying the Edge Embedding Lemma followed by the Subgroup Theorem.

Finally, let $n=kp+3$ for some $k$.  Embed $kp+1$ vertices as we did when $n=kp+1$.  Now add two additional vertices $v$ and $w$ in $C_g-C_f$ so that $f$ interchanges them.  Now $C_g-\{x,v,w\}$ has three components, whose closures will be the arcs $A_{xv}$, $A_{vw}$, and $A_{wx}$ required by hypothesis (3) of the Edge Embedding Lemma.  Now all of the hypotheses of the Edge Embedding Lemma are satisfied.  Hence we again obtain the required embedding of $K_n$ by applying the Edge Embedding Lemma followed by the Subgroup Theorem.

To get embeddings $\Lambda$ with $\TSG(\Lambda)=\Z_p$, we apply the Subgroup Theorem to each of the above embeddings.  

Finally, to get embeddings whose topological  symmetry group is $\Z_2$ let $p=3$.  Then for some integer $k$, $n=kp$, $kp+1$, or $kp+2$.  Hence by the above, $K_n$ has an embedding $\Gamma$ in $S^3$ with $\TSG(\Gamma)=D_p$.  Thus by the Subgroup Theorem, $K_{n}$ has an embedding $\Omega$ with $\TSG(\Omega)=\Z_2$.
\medskip

\noindent{\bf Case 2:}  $G=\Z_p\times \Z_q$ where $p$ and $q$ are odd and $p|q$.

In this case, let $G'\cong \Z_p\times \Z_q$ be a subgroup of $\Diff(S^3)$ generated by 
a rotation $g$ by $\frac{2\pi}{p}$ around a circle $C_g$ together with a rotation $f$ by $\frac{2\pi}{q}$ around a disjoint circle $C_f$ with $\mathrm{lk}(C_g, C_f)=1$.  Thus $C_f$ and $C_g$ are each setwise invariant under $G'$.  Also, the fixed point set of every non-trivial element of $G'$ is either the empty set, $C_g$, or $C_f$.

First let $k$ be an integer such that $n=kpq$.  Let $B$ denote a ball which is disjoint from $C_g\cup C_f$ and from its image under every non-trivial element of $G'$. Choose $k$ points in $B$.  The orbit of these points will be the $kpq$ embedded vertices of $K_n$.  Since no vertices are fixed or interchanged by any non-trivial element of $G'$, the hypotheses of the Edge Embedding Lemma are satisfied.  Thus by the Edge Embedding Lemma together with the Subgroup Theorem, we obtain the required embedding of $K_n$.

Next let $p=3$ and suppose that $n=kpq+3$ for some integer $k$.  We embed $kpq$ vertices asdescribed in the above paragraph.  Then add a vertex $v_1$ on $C_f$.  The orbit of $v_1$ under $G$ is 3 vertices $v_1$, $v_2$, and $v_3$ on $C_f$. Thus we have embedded all $kpq+3$ vertices.  Now it is not hard to check that the hypotheses of the Edge Embedding Lemma are satisfied.  Thus by the Edge Embedding Lemma together with the Subgroup Theorem, we obtain the required embedding of $K_n$.

Finally, let $p=q=3$ and suppose that $n=kpq+6$ for some integer $k$.  We embed $kpq+3$ vertices as in the above paragraph. Then add a vertex $w_1$ on $C_g$ together with its orbit.  This gives us three more vertices contained in $C_g$. Thus we have embedded all $kpq+6$ vertices.  Again it is not hard to check that the hypotheses of the Edge Embedding Lemma are satisfied. Now by the Edge Embedding Lemma together with the Subgroup Theorem, we obtain the required embedding of $K_n$.\end{proof}

\bigskip

\bigskip

\noindent {\bf Examples}

In order to demonstrate the usefulness of the above theorem, we apply the theorem to the complete graph $K_n$ for n=7, 15 and  27.

\begin{itemize}
\bigskip

\item $K_7$ can be embedded in $S^3$ as $\Gamma$ such that $G=\mathrm{TSG}_+(\Gamma)$ if and only if $G$ is one of the following:
$\Z_2$, $\Z_3$, $\Z_5$, $\Z_7$, $D_3$, $D_5$, or $D_7$.

\bigskip

\item $K_{15}$ can be embedded in $S^3$ as $\Gamma$ such that $G=\mathrm{TSG}_+(\Gamma)$ if and only if $G$ is one of the following:
$\Z_2$, $\Z_3$, $\Z_5$, $\Z_7$, $\Z_{13}$, $\Z_{15}$, $D_3$, $D_5$, $D_7$, $D_{13}$, $D_{15}$, or $\Z_3 \times \Z_3$.

\bigskip

\item $K_{27}$ can be embedded in $S^3$ as $\Gamma$ such that $G=\mathrm{TSG}_+(\Gamma)$ if and only if $G$ is one of the following:
$\Z_2$, $\Z_3$, $\Z_5$, $\Z_9$, $\Z_{13}$, $\Z_{25}$, $\Z_{27}$, $D_3$, $D_5$, $D_9$, $D_{13}$, $D_{25}$, $D_{27}$, $\Z_3 \times \Z_3$, or $\Z_3 \times \Z_9$.

\end{itemize}


\bigskip


\bigskip


\begin{thebibliography}{book}
\bibitem  {BLP}  M. Boileau, B. Leeb, J. Porti, \newblock {\it Geometrization of $3$-dimensional orbifolds}, \newblock Ann. of Math.  {\bf 162} (2005), 195--290.

\bibitem {Fl} E. Flapan,  {\it Rigidity of Graph Symmetries in the $3$-Sphere},  Journal of Knot Theory and its 	Ramifications, {\bf 4}, (1995), 373-388.

\bibitem{FMN1} E. Flapan, B. Mellor, R. Naimi, {\it Spatial Graphs with Local Knots}, preprint (2010).

\bibitem{FMN2} E. Flapan, B. Mellor, R. Naimi,  {\it Complete Graphs whose Topological Symmetry Groups are Polyhedral}, preprint (2010).

  \bibitem{FNPT}  E. Flapan, R. Naimi, J. Pommersheim, H. Tamvakis, {\it 
Topological Symmetry Groups of Embedded Graphs in the $3$-sphere}, 
Commentarii Mathematici Helvetici, {\bf 80}, (2005), 317-354.



\bibitem{FNT}
  	E. Flapan, R. Naimi, H. Tamvakis,  {\it 
Topological Symmetry Groups of Complete Graphs in the $3$-Sphere}, Journal of the
London Mathematical Society, {\bf 73}, (2006), 237-251.

\bibitem{Sm} P. A. Smith, {\it Transformations of finite period
II}, Annals of Math. {\bf 40} (1939), 690--711.

\end{thebibliography}
\end{document}